\documentclass[12pt,a4paper]{article}
\usepackage{graphicx}
\usepackage{amsmath}
\usepackage{amsfonts}
\newcommand{\be}{\begin{equation}}
\newcommand{\ee}{\end{equation}}
\newcommand{\ben}{\begin{enumerate}}
\newcommand{\een}{\end{enumerate}}
\newcommand{\beq}{\begin{eqnarray}}
\newcommand{\eeq}{\end{eqnarray}}
\newcommand{\beqn}{\begin{eqnarray*}}
\newcommand{\eeqn}{\end{eqnarray*}}

\newcommand{\bpf}{\begin{proof}}
\newcommand{\epf}{\end{proof}}
\newcommand{\bl}{\begin{lem}}
\newcommand{\el}{\end{lem}}
\newcommand{\bp}{\begin{prop}}
\newcommand{\ep}{\end{prop}}
\newcommand{\bd}{\begin{defn}}
\newcommand{\ed}{\end{defn}}
\newcommand{\bt}{\begin{thm}}
\newcommand{\et}{\end{thm}}
\newcommand{\R}{I\!\! R}
\def\nn{\nonumber}

\begin{document}
\title{On Sobolev spaces  and density theorems on Finsler manifolds}
\author{ Behroz Bidabad  and Alireza Shahi\thanks{Corresponding author} \\
{\small  Faculty of Mathematics, Amirkabir University of Technology,}\\
{\small 424  Hafez Ave. 15914, Tehran, Iran.}\\ {\small alirezashahi@aut.ac.ir;  bidabad@aut.ac.ir} }
\date{}
\maketitle
\begin{abstract}
Let $(M,F)$ be a $C^\infty$ Finsler manifold, $p\geq 1$  a real number, $k$  a positive integer and $H_k^p (M)$ a certain Sobolev space determined by a Finsler structure $F$.
 Here, it is shown that the set of all real $C^{\infty}$ functions  with compact support on $M$
is dense in the Sobolev space $H_1^p (M)$.
 This result  permits to approximate  certain solution of Dirichlet problem living on $H_1^p (M)$ by $C^ \infty$ functions with compact support on $(M,F)$.
 Moreover, let $W \subset M$ be a regular domain with  the $C^r$ boundary $\partial W$, then the set of all real functions in $C^r (W) \cap C^0 (\overline W)$ is dense in $H_k^p (W)$, where  $k\leq r$. This work is an extension of some density theorems of T. Aubin on Riemannian manifolds.
 \end{abstract}
Keywords:  Dirichlet problem, Density theorem, Sobolev spaces, Finsler manifolds.\\
MSC{   53C60,  46E35}
\section{Introduction}
\label{intro}
A Sobolev space is a vector space of functions endowed with a norm which is a combination of $L^p -$norm of the function itself and its derivatives up to a certain order.
Its objectives are to deal with some basic problems in geometry.
 A prominence application comes from the fact that solution of certain partial differential equations are naturally found in the Sobolev spaces, rather than in spaces of continuous functions. For instance solutions of the Helmholtz equation $\Delta u +k^2 u=0,$ are not smooth and are living on a (local) Sobolev space. As another application the Dirichlet problem consists of seeking a harmonic function in  Sobolev spaces. Indeed Sobolev spaces are devoted to  the study of properties of certain Banach spaces of weakly differentiable functions of several real variables which arise in connection with the numerous problems in theory of partial differential equations.
 As another application it is natural to wonder whether a function can be approximated by another one with better properties.
  Density problems permits to investigate  conditions under which a function on a Sobolev space can be approximated by smooth functions.

Historically, one the  significant  density theorems is proved by S. B. Myers \cite{My} in 1954 for compact Riemannian manifolds and then in 1959 by M. Nakai \cite{Na} for finite-dimensional Riemannian manifolds. Next in 1976 T. Aubin has investigated density theorems on Riemannian manifolds, cf. \cite{Au1} and \cite{Au2}.
 Recently, in 2009  an extension of Myers-Nakai theorem to infinite-dimensional, complete Riemannian manifolds is given in \cite{GJR1}.
 A similar result for the so-called finite dimensional Riemann-Finsler manifolds is given in \cite{GJR2} in 2010.
 Next in 2011 the Myers-Nakai theorem is extended to the Finsler manifolds of class $C^k$, where $k \in N\cup \{\infty \} $, cf., \cite{JS}. It is noteworthy to recall that if a compact Finsler manifold has non-zero constant sectional curvature then it reduces to a Riemannian manifold, cf., \cite{Ak,Bi} and hence density theorems hold well.

Previously, Y. Ge and Z. Shen \cite{GS} defined canonical energy functional on Sobolev spaces and investigated the eigenvalues and eigenfunctions related to this functional on Finsler manifolds.
 Moreover, Z. Shen in \cite{Sh} states Sobolev constants and compares these constants with first eigenvalue related to Finsler spaces.

 Let $(M,F)$ be an $n-$dimensional $C^\infty$ Finsler manifold, $p\geq 1$  a real number, $k$  a positive integer and $H_k^p (M)$ certain Sobolev space determined by the Finsler structure $F$.
  Denote by $D(M)$  the set of all real $C^\infty$ functions with compact support on $M$ and let $ \mathop {H^p_1}\limits^{\circ~} (M)$  be the closure of $D(M)$ in $H^p_1 (M).$

 Here, in this work we prove the following density theorems on Finsler manifolds.\\

 \noindent {\bf Theorem 1}. {\em ~Let $(M,F)$ be a forward geodesically complete, connected and $C^\infty,$ reversible Finsler manifold, then $\mathop{H^p_1}\limits^{\circ~} (M)=H_1^p (M)$}.\\

  \noindent {\bf Corollary 1}. {\em ~Let $(M,F)$ be a compact, connected, $C^\infty,$ reversible Finsler manifold and $f:M\longrightarrow{\R}$  a real function for which $\int _M fdv_F =0$, then the weak solution $u$ of the Dirichlet equation $\Delta u=f$ can be approximated by $C^\infty$ functions with compact support on $M$.}\\

Let  $W \subset M$ be a regular domain with $C^r$ boundary $\partial W$,  then $(\overline W,F)$ is called a Finsler manifold with $C^r$ boundary.\\

\noindent {\bf Theorem 2. }{\em Let $(\overline W,F)$ be a compact Finsler manifold  with $C^r$ boundary. Then $C^r (\overline W)$ is dense in $H_k^p (W)$, for $k \leq r.$}\\

As a consequence of density theorems in Finsler manifolds, we can approximate  solutions of partial differential equations on the Sobolev spaces determined by $F$, with $C^\infty$ or $C^r$ functions on $M$. Two examples are given in each case.
\section{Preliminaries}

Let M be an $n$-dimensional $C^\infty$ manifold, $TM$  the bundle of tangent spaces  $T_x M$ at $x\in M$.
 A (globally defined) Finsler structure on $M$ is a function $F:TM\longrightarrow [0,\infty)$ with the following properties :
(1) Regularity: $F:TM\longrightarrow [0,\infty)$ is $C^\infty$ on the entire slit tangent bundle $TM \backslash 0 $;
(2) Positive homogeneity: $F(x,\lambda y )=\lambda F(x,y) ~ \forall \lambda >0$;
(3) Strong convexity: The $n\times n$ Hessian matrix $(g_{ij}) =\frac {1}{2} ([F^2]_{y^i y^j})$ is positive-definite at every point of $TM\backslash 0$.
A Finsler structure $F$ is said to be reversible or absolutely homogeneous if $F(x,-y) =F(x,y)$ for all $(x,y) \in TM,$ cf., \cite{Ak} and \cite{BCS}.

Let $M$ be  connected and $ \sigma:[a,b]\longrightarrow M$ a piecewise $C^\infty$ curve with the velocity $\frac {d \sigma}{dt} =\frac {d \sigma ^i}{dt} \frac {\partial}{ \partial x^i} \in T_ {\sigma (t)} (M).$
Its integral length is defined by $L(\sigma)= \int _a^b F(\sigma , \frac {d \sigma}{dt}) dt.$ For $p$ and $x$ $\in M$ denote by $ \Gamma (p,x)$ the collection of all piecewise $C^\infty$  curves $\sigma:[a,b]\longrightarrow M$ with $\sigma (a)=p$ and $\sigma (b)=x$ and by the metric distance $d(p,x)$ from $p$ to $x$
 \be    \label{Eq;distance}
  d(p,x)=\inf\limits_{\Gamma (p,x)} L(\sigma).
 \ee
\noindent {\bf Lemma 2.1} \cite{BCS} {\em  ~Let (M,F) be a Finsler manifold. At any point $p \in M$, there exists a local coordinate system $\phi :\overline U\longrightarrow {\R}^n$ such that the closure of U is compact, $\phi (p)=0$ and $\phi$ maps $U$ diffeomorphically onto an open ball of ${I\!\!R}^n.$}\\
\noindent {\bf Lemma 2.2 }  \cite{BCS} {\em Let $(M,F)$ be a reversible Finsler manifold, then $d(p,x)=d(x,p)$ and $d$ is a metric on $(M,F)$}.\\
A Finsler manifold is said to be forward (resp. backward) geodesically complete  if every geodesic $\gamma (t),~ a \leq t <b,$ parameterized to have constant Finslerian speed, can be extended to a geodesic on $a \leq t<\infty.$ (resp. $-\infty <t\leq b$). If the Finsler structure $F$ is reversible, then $d$ is symmetric. In this case, forward completeness is equivalent to backward completeness.
Compact Finsler manifolds  at the same time are both forward and backward complete, whether $d$ is symmetric or not.
Fix a point $p,$ if $x$ is enough close to $p$, we can write $x=exp_p (v).$ The map $x=exp_p (v)$ is only $C^1$ at $p,$ although it is $C^ \infty$ away from $p,$ cf. \cite{BCS}.\\

\noindent {\bf Proposition 2.3 }  \cite{Ak,BCS}  {\em ~Let $(M,F)$ be a Finsler manifold, where $F$ is $C^\infty$ on $TM \backslash 0,$ and is positively (but perhaps not absolutely) homogeneous of degree one. For any fixed point $p \in M,$ the function $d_p ^2$ is :\\
(a) $C^\infty$ in a punctured neighborhood of $p.$\\
(b) $C^1$ at $p$ and has zero derivative there.\\
(c) $C^2$ at $p$ if and only if $F$ is Riemannian.}

 \section{A Laplacian on Finsler manifolds and Osculating Riemannian metric}
 Let $(M,F)$ be a Finsler manifold. In this section, Busemann volume form defined in \cite{Bu} and  leading term of the Laplacian given in  \cite{Ce2} are used to determine an appropriate Riemannian metric tensor and a volume form on Finsler manifolds. At a point $x\in M$, let $\{b_{i}\}_{i=1}^{n}$ be a basis for $T_x M$ and $\{\theta^i\}_ {i=1}^ {n}$
 its dual basis. H. Busemann exhibited the following natural Finsler volume form on a Finsler space such that the measure defined by this volume form is just the Hausdorff measure of the induced distance function $d_F$, providing that $F$ is reversible.
 Let $dv_F =\sigma_F (x)\theta^1 \wedge \theta^2 \cdots \wedge \theta^n,$
 where $\sigma_F (x)=\frac {vol(B^n (1))} {vol((y^i)\in {\R}^n :F(x,y^i b_i )<1)}~,$ and $vol(B^n (1))$ denotes the Euclidean volume of unit ball.
  Let $U$ be an open subset of $M$ the volume of $U$ with respect to $F$ is defined by $vol_F (U)=\int_Udv_F$.
 A linear differential operator $A(u)$ of order $2m$ on $M$, written in a local chart $(V,\psi)$, is an expression of the form $A(u)=\sum _{l=0}^{2m} a_l ^{\alpha _1 \alpha _2 \cdots \alpha _l} \nabla _{\alpha _1 \alpha _2 \cdots \alpha _l} u$, where $a_l$ are $l-$tensors and $u \in C^{2m} (M).$ The adjoint of $A$ is defined by $A^* (\phi)=(-1)^l \nabla ^l (\phi a_l).$  We say that $u$ is a \textit{ weak solution} of $A(u)=f$ if for all $\phi \in D(M)$ we have $\int _M {uA^* (\phi)dv_F}=\int _M {f \phi dv_F}.$\\
  Assuming $F=\sqrt{g_{ij}y^i y^j}$ is Riemannian, $dv_F$ reduces to the common Riemannian volume form, $dv_\alpha =\sqrt{detg_{ij}} dx^1\wedge dx^2 \wedge \cdots\wedge dx^n,$ cf., \cite{CS}.
 A Laplacian for Finsler manifolds is introduced in \cite{Ce1} which leads to the definition of Osculating Riemannian metric.\\
\noindent {\bf  Theorem 3.1}  \cite{Ce1}{\em ~ Let $(M,F)$ be a Finsler manifold, $p \in M$  and $f:M\longrightarrow {\R}$  a smooth function, then the Laplacian of $f$ in $p$, exists and is a linear $2nd$ order differential operator with the following coordinate expression:
\be    \label{Eq;Laplace}
  \Delta(f)=(n+2)(\frac{\int_I y^i y^j dy}{\int_I dy}f_{ij} + lower~ order ~terms),
 \ee
 where $y=y^i \frac {\partial}{\partial x^i}\in T_x M$ , $I=\{y\in T_x M : F(y)\leq1\}$ and $f_{ij} =\frac{\partial ^2 f}{\partial x^i \partial x^j}.\\$}
   Coefficients of first term in the right hand side of (\ref{Eq;Laplace}) are symmetric, positive-definite, twice-contravariant tensor. Denoting by
 $$K^{ij}=(n+2)(\frac{\int_I y^i y^j dy}{\int_I dy}),$$
one can easily check that the coefficients $K^{ij}$ can be considered as the inverse component of a Riemannian metric called  Osculating Riemannian metric, cf., \cite{Ce3}.\\
Let
 $\alpha=\alpha _{i_1 \cdots i_p }dx^{i_1}\wedge dx^{i_2 }\wedge\cdots \wedge dx^{i_P }$ be a $p-$form on $M$ then $$|\alpha|^2 =\alpha _{i_1 \cdots i_p }\alpha_{j_1 \cdots j_p}K^{i_1 j_1}\cdots K^{i_p j_p }.$$
 For a real function $\phi\in C^r (M)$ where, $r\geq 0$ is an integer and $\nabla$  the  Levi-Civita connection on the associated Riemannian manifold $(M,K)$, we denote the covariant derivative of order $l$ by $\nabla ^l$ and consider $\nabla ^l\phi $ as a $l-$form on $M$. Thus, for $l\leq r$ one can define
$$|\nabla ^l \phi|^2 =K^{i_1 j_1 }\cdots K^{i_l j_l }\nabla_{i_1} \cdots \nabla _{i_l }\phi \nabla_{j_1}\cdots \nabla _{j_l} \phi.$$
For instance $|\nabla^1 \phi |^2 =K^{ij} \nabla_i \phi \nabla_j \phi,$
where $\nabla \phi$ is the first covariant derivative of $\phi.$
\section{A family of Sobolev spaces on Finsler manifolds}
In this section  we use the Yang's method to define certain Sobolev spaces on Finsler manifolds. Let $dv_F$ be the Busemann volume form,  $f \in C^\infty (M)$ and denote
$$L^p (M):=\{f:M\longrightarrow {\R}~ is~measurable~ |~ \int _M f^p (x)dv_F <\infty \},$$
and $\parallel f\parallel_{L^p}: = (\int_M |f|^p  dv_F )^{\frac {1}{p}}.$\\
Let us consider the vector space $C_k ^p$ of $C^\infty $ functions $\phi$, such that $|\nabla ^l \phi |\in L^p (M)$, for all integers $k$ and $l$  with $0\leq l\leq k$, and real number $p\geq 1$.
The Sobolev space $H_k ^p (M)$ is defined to be a completion of $C_k ^p $ with respect to the norm $\parallel \phi \parallel _{H_k ^p (M)} :=\sum _{l=o} ^k \parallel  \nabla ^l \phi \parallel _p .$
 $\mathop{H^p_1}\limits^{\circ~} (M)$ is the closure of $D(M)$ in $H_1 ^p (M)$, where $D(M)$ is the space of $C^\infty $  functions with compact support in $M$ and $H_0 ^p =L^p.$\\
T. Aubin has extended the notion of Euclidean Sobolev spaces to Riemannian manifolds. Y. Yang  generalized this definition to Finsler manifolds using Busemann volume form and Osculating Riemannian metric cf., \cite{Ya}. We use this definition of Sobolev spaces in Finsler manifolds to prove two density theorems.\\
The Sobolev's norm may be also defined to be $(\sum_{l=o}^k \parallel\nabla^l \phi \parallel  _p ^p )^\frac {1}{p}$, which is equivalent to the norm used in this paper defined by $\parallel \phi \parallel _{H_k ^p (M)}$.

Let $J$ be a nonnegative, real-valued function, on the space of $C^\infty $ functions with compact support on ${\R}^n,$ denoted by $C_0^\infty ({\R}^n)$ and having properties :\\
(i)$ J(x)=0$ if $| x| \geq 1 .$\\
(ii)$ \int  _{{\R}^n} J(x)dx =1.$\\
Consider the function $J_\epsilon (x)=\epsilon ^{-n} J(\frac{x}{\epsilon})$ which is nonnegative in  $C_0^\infty ({\R}^n)$ and satisfies :\\
(i) $J_\epsilon (x) =0$ if $| x| \geq \epsilon ,$\\
(ii) $ \int  _{{\R}^n} J_\epsilon (x)dx =1.$\\
$J_\epsilon$ is called a mollifier and the convolution $J_\epsilon * u(x) := \int _{{\R}^n} J_\epsilon (x-y) u(y)dy$ defined for the function $u$ is called a mullification or regularization of $u.$\\
\noindent {\bf Lemma 4.1 }  \cite{Ad} {\em ~  Let u be a function defined on $\Omega \subset {\R}^n$ and vanishes identically outside the domain $\Omega :$\\
(a) If $u \in L^1_{Loc}( \overline {\Omega }) $ then $J_\epsilon * u \in C_0^\infty ({\R}^n).$\\
(b) If also $supp(u)\subset  \Omega,$ then $J_\epsilon * u \in C^\infty _0 (\Omega)$ provided $\epsilon < dist(supp(u),\partial \Omega).$\\
(c) If $u \in L^p$ where $1\leq p< \infty$ then $J_\epsilon * u \in L^p (\Omega ).$ Moreover,\\ $\parallel J_\epsilon *u  \parallel _p \leq \parallel u \parallel _p$
and $\lim\limits _{\epsilon \longrightarrow 0^+} \parallel J_\epsilon *u -u \parallel _p =0$\\}

\noindent {\bf Lemma 4.2 } \cite{Ad} {\em ~ Let $u \in C_k^p (\Omega)$ where $1\leq p < \infty $ and $k$ be an integer. If $\Omega' \subset \Omega,$ then $\lim\limits _{\epsilon\longrightarrow o^+} J_\epsilon * u(x)=u(x)$ in $C^p_k (\Omega' ).$}\\

\noindent {\bf Lemma 4.3 } \cite{Ad} {\em ~ Let  $u \in H_k^p (\Omega)$ and $1 \leq p<\infty.$ If $\Omega' \subset  \Omega,$ then $\lim\limits_{\epsilon\longrightarrow o^+} J_\epsilon * u=u,$ in $H_k^p (\Omega')$}
\section{Proof of the density theorems on Finsler manifolds}
Let $(M,F)$ be a Finsler manifold and $D(M)$ the space of $C^\infty$ functions with compact support on $M.$
  In this section  we use Hopf-Rinow's theorem to introduce  the first density theorem on boundary-less Finsler manifolds and investigate another density theorem for Finslerian manifolds with $C^r$ boundary.
\subsection{Case of manifolds without boundary}
\noindent {\bf  Proof of Theorem 1.}  We first recall that, by  definition $\mathop{H^p_1}\limits^{\circ~} (M) $ is the closure of $D(M)$ and $C^\infty (M)\cap  H_1^p (M)$ is dense in $H^p_1 (M).$ To prove this Theorem  we have  to show that for each $C^\infty$ real function $\phi$ on $M$, $\forall x\in M,$ $\phi (x)$ is in $H^p _1 (M)$, $\phi (x)$ can be approximated by $C^\infty$ functions with compact support in $M.$\\
 Let $\phi\in C^\infty (M) \cap H_1^p (M)$ and  fix the point $x_0\in M.$ By Hopf-Rinow's~theorem on forward geodesically complete Finsler manifolds, every pair of points containing $x_0$ can be joined by a minimizing geodesic emanating from $x_0$. Consider the Finslerian distance function defined by (\ref{Eq;distance}), from $x_0$ and sequence of the functions $\phi_j (x) =\phi (x) f(d(x_0 ,x)-j),$ where $f:{\R}\longrightarrow {\R} $ is defined by
 \begin{equation}
    f(t)=\left\{ \begin{array}{ll}1& t\leq 0\\ 1-t&0<t<1\\ 0&t\geq 1\end{array}.\right.
 \end{equation}
$f$ is continuous decreasing  function on ${\R}$ and  differentiable almost everywhere. Next we should prove, each $\phi _j (x)$ has a compact support.\\
If $d(x_0 ,x) \leq j$ then $d(x_o ,x)-j\leq 0$ or $f(d(x_o ,x)-j)=1$ hence $\phi _j (x)=\phi (x).$
If $d(x_0 ,x)\geq j+1$ then $d(x_o ,x)-j\geq 1$ or $\phi _j (x)=0.$\\
 Thus each $\phi _j$ has a compact support.
$\phi (x)$ is $C^ \infty $ real function on $M$ and $\phi (x)\in H^p_1 (M)$, therefore $\nabla \phi (x)$ exists and is bounded almost everywhere on $M$.\\
It is well known, the Finslerian distance function $d$ is $C^\infty$ out of a small neighborhood of $x_0$ and  $C^1$ in a punctured neighborhood of $x_0$ cf., \cite{Wi}.
Therefore $\phi _j (x)$ are Lipschitzian  on $M$.
We have $\lim\limits _{\scriptscriptstyle j\longrightarrow \infty} {\phi _j (x)}=\phi (x)$
and by definition  of $\phi _j$, $|\phi _j (x)| \leq |\phi (x)|$ where $\phi (x)$ is $C^\infty$ function which lies in $L^p$. Therefore by the Lebesgue dominated convergence theorem $\phi _j (x)$ and $\phi _j -\phi $   are also in $L^p $. $(M,F)$ is reversible therefore $d$ is symmetric. Consider forward metric ball
$B_{x_0}^+ (j)=\{x \in M: d(x_0 ,x)<j\}$, we have\\
 $(\int_ M |\phi _j -\phi|^p dv_F)^\frac{1}{p} \leq(\int _{B_{x_0}^+ (j)} |\phi _j -\phi|^p dv_F)^\frac{1}{p} +(\int_{M \setminus{{B_{x_0}^+} (j)}} |\phi _j -\phi|^p dv_F)^\frac{1}{p}. $\\
 On $B_{x_0}^+ (j)$ we have $d(x_0 ,x)<j$ or $\phi_j (x)=\phi (x)$ and when $j$ tends to infinity $M\setminus {B_{x_0}^+ (j)}=\emptyset$ therefore

$(\int_ M |\phi _j -\phi|^p dv_F)^\frac{1}{p} \leq (\int_{M\setminus {B_{x_0}^+ (j)}} (2|\phi |)^p dv_F)^\frac {1}{p} =2 (\int_{M\setminus {B_{x_0}^+ (j)}} (|\phi |)^p dv_F)^\frac {1}{p}$, hence
$$\parallel \phi _j - \phi \parallel _p =(\int _M|\phi_j -\phi|^p dv_F)^{\frac{1}{p}} \longrightarrow 0.$$
Next we prove $\parallel \phi _j - \phi \parallel _{H^p_1 (M)} $ converges to zero. To this end it suffices to show that $\parallel \nabla \phi _j -\nabla \phi \parallel _p $ or $ \parallel \nabla (\phi_j - \phi) \parallel _p $ converges to zero .
We use Leibnitz's formula, $\phi _j (x)=\phi (x) f(d(x_0,x)-j)$, which leads to
$$| \nabla \phi _j (x)| \leq | \nabla \phi (x) | + | \phi (x) | \sup\limits _{t \in [0,1]} | f'(t) |.$$
Again with Lebesgue dominated theorem, we have  $|\nabla \phi _j| \in L^p$ and hence $\phi _j (x) \in H_1^p (M).$
\begin{equation}\nn
(\int_ M |\nabla \phi _j -\nabla \phi|^p dv_F)^\frac{1}{p} \leq (\int_ {B_{x_0}^+ (j)} |\nabla\phi _j -\nabla\phi|^p dv_F)^\frac{1}{p}+ (\int_ {M\setminus {B_{x_0}^+ (j)}} |\nabla \phi _j -\nabla \phi|^p dv_F)^\frac{1}{p},
 \end{equation}
 by definition of $\phi_j$ and $f(t)$ and Leibnitz's formula we have
 \begin{eqnarray}\nn
(\int_ M |\nabla \phi _j -\nabla \phi|^p dv_F)^\frac{1}{p}&\leq &(\int_ {M\setminus {B_{x_0}^+ (j)}} |\nabla \phi _j -\nabla \phi|^p dv_F)^\frac{1}{p}\\
&\leq&(\int_ {M\setminus {B_{x_0}^+ (j)}} |\nabla \phi |^p dv_F)^\frac{1}{p} + (\int_ {M\setminus {B_{x_0}^+ (j)}} | \phi|^p dv_F)^\frac{1}{p},\nn
\end{eqnarray}
hence
$$\parallel \nabla \phi _j-\nabla \phi \parallel _p =\parallel \nabla (\phi_j - \phi) \parallel _p = (\int _M|\nabla  (\phi _j -\phi)|^p dv_F )^{\frac {1}{p}} \longrightarrow 0.$$
Therefore
 $$\parallel \phi _j -\phi \parallel _ {H_1^p (M)}= \parallel \phi _j - \phi \parallel _p + \parallel \nabla(\phi _j -\phi )\parallel _p \longrightarrow 0,$$ thus $\phi_j$ converges to $\phi$ in $H_1^p (M)$.
Now it remains to approximate each $\phi _j$ by functions in $D(M)$. Let $j$ be a fixed index, for which $\phi _j$ has a compact support. Let $K$ be the compact support of  $\phi _j$ and  $\{V_i \}_{i=1}^m $ be a finite covering of $K$ such that by means of Lemma 2.1, for fix index $i$, $V_i$ is homeomorphic to the open unit ball $B$ of ${\R}^n$ . Let $(V_i ,\psi _i )$ be the corresponding chart, we complete the proof by partition of unity's theorem. Let $\{\alpha _i\}$ be a partition of unity subordinate to the covering $\{V_i \}_{i=1}^m $. Thus for approximating $\phi _j$ by $C^ \infty$ functions  with compact support in $H_1^p (M)$, it remains to approximate each $\alpha _i \phi _j$ for $1 \leq i\leq m.$ For fixed $i,$  $\psi _i$ is a homeomorphic map between $V_i$ and the unit ball $B$. Consider the functions $(\alpha _i \phi _j)\circ \psi_i ^{-1}$ which have their supports in $B$. Let we denote $(\alpha _i \phi _j)\circ \psi_i ^{-1}$ by $u$, to be consistent with the notation of Lemmas 4.1, 4.2 and 4.3. \\ Consider the convolution $J_\epsilon *u$ with $\lim\limits _{\epsilon \longrightarrow o^+} J_\epsilon *u =u.$ Let $h'_\epsilon =J_\epsilon *u \in C^\infty ({\R}^n),$ then $h'_\epsilon$ has a compact support, that is  $h'_\epsilon \in D(B)$. We approximate $u=(\alpha _i \phi _j)\circ \psi_i ^{-1}$ by $h_k =h'_\frac {\epsilon }{k}$. More precisely
$$\lim\limits _{k\longrightarrow \infty } h_k = \lim\limits_{k\longrightarrow \infty} h'_\frac{\epsilon}{k}=\lim\limits_{\epsilon \longrightarrow o^+} h'_\epsilon=u= (\alpha _i \phi _j)\circ \psi_i ^{-1}.$$ Also $h_k$ is $C^\infty$, $u$ is Lipschitzian  and $\nabla h_k \longrightarrow \nabla u.$
Hence $h_k$ converges to $h$ in $H_1^p (B).$\\
Now $h_k \circ\psi _i$ converges to $\alpha _i \phi _j$ in $H_1^p (V_i)$ and $h_k \circ\psi _i \in D(V_i ).$
Thus we have approximated each $\alpha _i \phi _j$ by functions in $D(V_i )$. This completes the proof.
\hspace{\stretch{1}}$\Box$\\

It should be remarked that theorem 1 is also true for compact, connected and reversible Finsler manifolds. In fact compact Finsler manifolds are forward geodesically complete.\\

\noindent {\bf  Proof of Corollary 1.} Given a compact Finsler manifold, using the Aubin's technic on Riemannian manifolds, Yang on a non-published work \cite{Ya}  has proved  the Dirichlet equation $\Delta u=f$ has a weak solution $u$ on $H_1^2 (M)$, providing that $\int _M fdv_F =0$ . Theorem 1 shows that this solution can be approximated by  $C^ \infty$ functions with compact support on $M$. This completes the proof.\hspace{\stretch{1}} $\Box$\\

 \noindent {\bf Example 5.1.1 }   Let $M=S^2$ and $F$ be a reversible Finsler structure on $S^2$. Define $f: S^2\rightarrow {\R}$  by
  \begin{equation}
    f(x_1 ,x_2,x_3)=\left\{\begin{array}{ll}1&  \quad x_3> 0\\-1& \quad x_3\leq 0\end{array}.\right.
 \end{equation}
We can see easily $\int_M fdv_F=0$. Hence the Dirichlet problem $\Delta u =f$ has a weak solution $u \in H_1^2 (S^2)$. $S^2$ is compact, connected  and $F$ is reversible, therefore by means of Corollary 1,  we can approximate the weak solutions of $\Delta u =f$ by $C^\infty$ functions with compact support on $S^2.$\\
 \subsection{Case of manifolds with $C^r$ boundary}
In proof of Theorem 2 we use the technic applied in proof of the following theorem on half-spaces on ${\R}^n.$\\
\noindent {\bf Theorem 5.1 }  \cite{Au2}{\em ~ $C^\infty (\overline E)$ is dense in $H^p_k (E)$, where E is a half-space  $E=\{x \in {\R}^n : x_1 <0\}$ and $C^\infty (\overline E)$ is the set of functions that are restriction to $\overline E$ of $C^\infty$ functions on ${\R}^n.$}\\

\noindent {\bf Proof of Theorem 2. } Let $\phi$ be a real $C^\infty$ function on the Sobolev space  $ H_k^p (W)$, that is, $\phi \in C^ \infty (W) \cap H_k^p (W). $
Here we approximate $\phi$ by the functions in $C^r (\overline W)$. Since $\overline W$ is compact, we can consider $(V_i,\psi _i ),i=1,\cdots ,N$ as a finite $C^r$ atlas on $\overline W.$ Each $V_i$ depending on  $V_i \subset W$ or  $V_i$ has intersection with $\partial W,$ is homeomorphic either to the unit ball of ${\R}^n$ or  a half-ball $D=B\cap \overline E,$ where $E=\{x\in {\R}^n : x_1 <0\}.$  Here $W$ is not complete, but $\overline W$ is compact, so we can use the partition of unity theorem. Let $\{\alpha _i\}$ be a $C^\infty$ partition of unity subordinate to the finite covering $\{V_i\}_{i=1}^m$ of $\overline W.$ By partition of unity's properties it remains to show that each $\alpha _i \phi$,  supported in $V_i,$ can be approximated by functions in $C^r (V_i).$ Each $V_i$ is homeomorphic either to the unit ball $B$ or a half-ball $D.$ First let  $V_i$ be homeomorphic to $B$ and let $\psi _i$ be a homeomorphism  between $B$ and $V_i.$ Consider $(\alpha _i \phi)\circ\psi_i ^{-1} :=h,$ the support of $h$ is in $B.$ Similar to the same argument in the proof of Theorem 1, for derivative of higher orders we use Leibniz's formula to show that $\parallel h_l \circ\psi _i - \alpha _i \phi \parallel _ {H_k^p (V_i)}$ converges to zero.

By appropriate choice of $V_i,$ $ h_l \circ\psi _i$ and $\alpha _i \phi$ which are $C^\infty$ real functions on $V_i,$  we can use covariant derivative of these functions up to order $\leq r$ on $(W,F)$ and theorem holds well in this case. Hence $\parallel h_l \circ\psi _i - \alpha _i \phi \parallel _ {H_k^p (V_i)}$ converges to zero.\\
Now let $V_i$ be homeomorphic to $D=B\cap \overline E$ and $\psi _i$ be a homeomorphism  between $V_i$ and $D.$ Consider the sequence of functions  restricted to $D$ where $((\alpha _i \phi )\circ\psi_i^{-1})(x_1 -\frac{1}{m},x_2 ,x_3, \cdots ,x_m )$ converges to $((\alpha _i \phi )\circ\psi_i^{-1})$ in $H^p_k (D).$ Let $\phi \in H_k^p (W) \cap C^\infty (W),$ by appropriate choice of $\{V_i\}$ and $\{\alpha _i\}, $ the restriction of $((\alpha _i \phi )\circ\psi_i^{-1})(x_1 -\frac{1}{m},x_2 ,x_3, \cdots ,x_m )$ to $D$ and its derivative up to order $r$  converge to $((\alpha _i \phi )\circ\psi_i^{-1})$ in $H_k^p (D)$, where  $D$ has the Euclidean metric. In this case the Osculating Riemannian metric and all its derivative of order $\leq r$ are bounded on $V_i$ . Indeed
\begin{eqnarray*}
| \nabla ^s \phi | ^p& =&(| \nabla ^s \phi | ^2 )^ \frac {p}{2}\\
&=&(K^{i_1 j_1}\cdots K^{i_s j_s}\nabla _{i_1}\cdots \nabla _{i_s} \phi \nabla _{j_1}\cdots \nabla _{j_s} \phi )^ \frac {p}{2},
\end{eqnarray*}
is bounded on $V_i$, $1 \leq s \leq  r,$ and
\begin{eqnarray*}
&&\int _{V_i}|\nabla^s(h_l\circ\psi_i-\alpha_i\phi)|^pdv_F
\\
&&=\int_{V_i}(K^{i_1j_1}\cdots K^{i_sj_s}\nabla_{i_1}\cdots\nabla_{i_s}(h_l\circ\psi_i-\alpha_i\phi)\nabla_{j_1}
\cdots\nabla_{j_s}(h_l\circ\psi_i-\alpha_i\phi))^\frac{p}{2}dv_F\longrightarrow 0,
 \end{eqnarray*}
 Hence  $h_l \circ\psi _i \longrightarrow \alpha _i \phi $ in $H_k^p (V_i),~k\leq r$, $h_l \circ\psi _i \in C^r (\overline W)$  and proof is complete.~  
\hspace{\stretch{1}}$\Box$\\

In the following example we can see that the assumption $k\leq r$ in above theorem can not be omitted.\\
\noindent {\bf Example 5.2.1 } Let $\overline W=[-1,1]\times[0,1]$ be a manifold with boundary of class $C^0$. Denote the points of $W$ by $x=(x^1 ,x^2)$ and the points of $T_x (W)$ by $y=(y^1 , y^2)$. Let $F(x,y)= \sqrt{g_x (y,y)}$ be a Finsler structure defined by
$g=g_{ij} ~dx^i \otimes dx^j=(dx^1)^2 + (dx^2)^2$ on $W$ and $dv_F$ be the Busemann volume form given by $dv_F =\sigma _F (x)~ dx^1 \wedge dx^2 =dx^1 \wedge dx^2$, where $\sigma_F(x)=\sqrt{\det g_{ij}}$.  Define $u:W\rightarrow {\R}$ by
 \begin{equation}
 u(x^1 ,x^2)=\left\{ \begin{array}{ll}1& \quad x^1> 0 \\ 0& \quad x^1\leq 0\end{array}.\right.
 \end{equation}
 $u$ belongs to $H_1^p (W).$ We claim that for sufficiently small $\varepsilon$, there is no $\phi \in C^1 (\overline W)$ such that $\parallel u-\phi \parallel_{H^p_1 (W)} < \varepsilon.$ Assume for a while that our assumption is not true and the function $\phi$ exists.  Let $S=\{(x^1 , x^2) : -1 \leq x^1 \leq 0 ~ , 0 \leq x^2 \leq 1\}$  and $K= \{(x^1 , x^2) : 0<x^1 \leq 1 ~ , 0<x^2 \leq 1\}$, then $\overline W = S\cup K.$ On $S$, $u(x^1,x^2)=0$, hence $\parallel 0-\phi \parallel _{H_1^p (S)} < \varepsilon $ or $\parallel \phi \parallel _{L^1 (S)} + \parallel \nabla\phi \parallel _{L^1 (S)} < \varepsilon ,$ therefore $\parallel \phi \parallel _{L^1 (S)} < \varepsilon$. On $K$, $u(x^1 ,x^2)=1$ thus $\parallel 1-\phi \parallel_{L^1 (K)} < \varepsilon$ or $\parallel \phi \parallel_{L^1 (K)} > 1-\varepsilon$. Put $\psi (x^1)= \int _0 ^1 \phi (x^1 , x^2)dx^2$, then there exist $a$ and $b$ with $-1 \leq a \leq 0$ and $0<b \leq 1$ such that $\psi (a)<\varepsilon$ and $\psi (b)>1-\varepsilon$. Thus $$1-2\varepsilon < \psi (b)- \psi (a) = \int _a ^b \psi' (x^1)dx^1 \leq \int _{\overline W} |D_{x^1} \phi (x^1 , x^2)| dx^1 dx^2 \leq$$
 $$ \big(\int _{\overline W} 1^{p'} dx^1 dx^2 \big)^ \frac {1}{p'} \big(\int _{\overline W} |D_{x^1} \phi (x^1 , x^2)|^p dx^1 dx^2 \big)^ \frac {1}{p} =2^{\frac{1}{p'}} \parallel D_{x^1} \phi (x^1 , x^2) \parallel _{L^p (\overline W)} $$$$<2^{\frac{1}{p'}} \varepsilon,$$ where $\frac{1}{p} + \frac{1}{p'} =1$. Hence $1<(2+2^{\frac {1}{p'}})\varepsilon$ which is not possible for small $\varepsilon$. This contracdict  our provisional assumption and prove our statement.


%




\end{document}